\documentclass[11pt,leqno,fleqn]{article}
\usepackage{amsmath,amssymb}
\usepackage{psfrag}
\newcounter{sectie}
\newcommand{\sectz}[1]{\refstepcounter{sectie}
\section*{\boldmath \thesectie. #1}%
}
\newcommand{\dyz}[1]{%
\refstepcounter{equation}%
\begin{list}{}{
\topsep 3mm
\leftmargin 18mm
\rightmargin 0cm
\itemsep 0mm
\listparindent 0mm
\parsep 0mm
\itemsep 0mm
\labelsep 0mm
\labelwidth 18mm
}%
\item[\rm (\theequation)\hfill]
#1
\end{list}%
}
\newcommand{\bgcirc}{\raisebox{1.4pt}{\mbox{\scriptsize $\mathbf{\bigcirc}$}}}
\newcommand{\GG}{{\cal G}}
\newcommand{\oR}{{\mathbb{R}}}
\newcommand{\oZ}{{\mathbb{Z}}}
\newcommand{\dez}[1]{\dyz{\raggedright$\displaystyle #1 $}}
\newcommand{\RR}{{\cal R}}
\newcommand{\de}[2]{\dy{#1}{\raggedright$\displaystyle #2 $}}
\newcommand{\dy}[2]{%
\refstepcounter{equation}%
\label{#1}%
\begin{list}{}{
\topsep 3mm
\leftmargin 18mm
\rightmargin 0cm
\itemsep 0mm
\listparindent 0mm
\parsep 0mm
\itemsep 0mm
\labelsep 0mm
\labelwidth 18mm
}%
\item[\rm (\theequation)\hfill]
#2
\end{list}%
}
\newcommand{\di}[2]{%
\refstepcounter{equation}%
\label{#1}%
\begin{list}{}{
\topsep 5mm
\leftmargin 10mm
\rightmargin 0cm
\itemsep 0mm
\listparindent 0mm
\parsep 0mm
\labelsep 1mm
\labelwidth 10mm
}%
\item[\rm (\theequation)\hfill]
\begin{list}{}{
\topsep 0mm
\leftmargin 8mm
\rightmargin 0mm
\itemsep 0mm
\listparindent 0mm
\parsep 0mm
\labelsep 1.5mm
\labelwidth 6.5mm
}
#2
\end{list}%
\end{list}%
}
\newcommand{\nr}[1]{\item[{\rm (#1)}]}
\newcommand{\nrs}[1]{\item[{\rm (#1)}]\vspace{-\itemsep}}
\newcommand{\dps}{\displaystyle}
\newcommand{\rf}[1]{{\rm (\ref{#1})}}
\renewcommand{\phi}{\varphi}
\newcommand{\join}[1]{\raisebox{-.05\height}{\mbox{\hspace*{2pt}\footnotesize$\overset{\hspace*{0.5pt}{#1}}{{\scriptsize\vee}}$\hspace*{2pt}}}}
\newcommand{\dyyz}[1]{\dyz{\raggedright$\dps#1$}}
\newcommand{\thmnn}[1]{\vspace{4mm}\noindent{\bf Theorem.}{\it #1}}
\newcommand{\OO}{{\cal O}}
\newcommand{\lemmann}[1]{\vspace{4mm}\noindent{\bf Lemma.}{\it #1}}
\newcommand{\pf}{\vspace{3mm}\noindent{\bf Proof.}\ }
\newcommand{\MM}{{\cal M}}
\newcommand{\id}{\text{\rm id}}
\newcommand{\T}{^{\sf T}}
\newcommand{\sgn}{\text{\rm sgn}}
\newcommand{\minus}{{\hspace*{-1pt}-\hspace*{-1pt}}}
\newcommand{\bx}{\hspace*{\fill} \hbox{\hskip 1pt \vrule width 4pt height 8pt depth 1.5pt \hskip 1pt}

\addvspace{4mm}}
\newcommand{\dyy}[2]{\dy{#1}{\raggedright$\dps#2$}}
\begin{document}

\begin{center}
{\large\bf ON THE EXISTENCE OF REAL R-MATRICES FOR VIRTUAL LINK INVARIANTS

}
\vspace{4mm}
Guus Regts\footnote{ University of Amsterdam.
The research leading to these results has received funding from the European Research Council
under the European Union's Seventh Framework Programme (FP7/2007-2013) / ERC grant agreement
n$\mbox{}^{\circ}$ 339109.},
Alexander Schrijver$\mbox{}^1$,
Bart Sevenster$\mbox{}^1$

\end{center}

\noindent
{\small
{\bf Abstract.}
We characterize the virtual link invariants that can be described as partition function of a real-valued R-matrix,
by being weakly reflection positive.
Weak reflection positivity is defined in terms of joining virtual link diagrams, which is a specialization
of joining virtual link diagram tangles.
Basic techniques are the first fundamental theorem of invariant theory, the Hanlon-Wales theorem on
the decomposition of Brauer algebras, and the Procesi-Schwarz theorem on inequalities for closed orbits.

}

\sectz{Introduction}

This paper is inspired by some recent results in the range of
characterizing combinatorial parameters using invariant theory, in particular by Szegedy [12] and
Freedman, Lov\'asz, and Schrijver [1].
We here consider the application to virtual links, which requires some new techniques from the representation
theory of the symmetric group.
The concepts of virtual link diagram and virtual link were introduced
by Kauffman [5]; see Manturov and Ilyutko [7] and Kauffman [6] for more background.

A {\em virtual link diagram} is an undirected 4-regular graph $G$ such that
at each vertex $v$ a cyclic order of the edges incident with $v$ is specified, together
with one pair of edges opposite at $v$ that is labeled as `overcrossing'.
The standard way of indicating this is as
\dyz{
\raisebox{-.35\height}{\scalebox{0.18}{\includegraphics{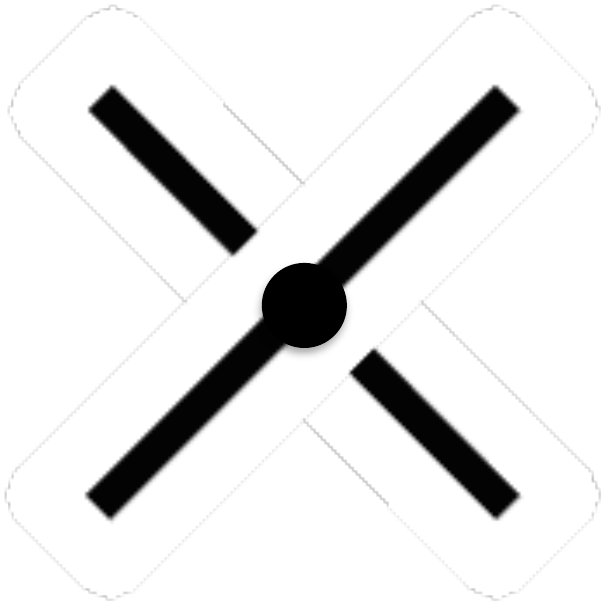}}}
~~~
or just
~~~
\raisebox{-.35\height}{\scalebox{0.18}{\includegraphics{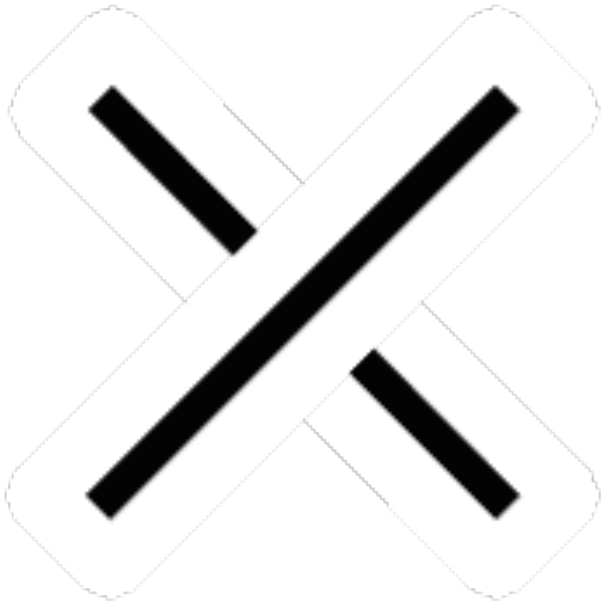}}}
.
}
Vertices of a virtual link diagram are called {\em crossings}.
Loops and multiple edges are allowed.
Moreover, the `unknot' is allowed, that is, the loop $\bgcirc$ without a crossing.
Let $\GG$ denote the collection of virtual link diagrams,
two of them being the same if they are isomorphic.

In the usual way, Reidemeister moves yield an equivalence relation on virtual link diagrams.
A {\em virtual link} is an equivalence class of virtual link diagrams.
A {\em virtual link invariant} is a function defined on $\GG$
that is invariant under Reidemeister moves.
(So in fact it is a function on virtual links, but the
definition as given turns out to be more convenient.)

A virtual link diagram can be seen as the projection
of a link in $M\times\oR$ on $M$, where $M$ is some oriented surface.
Since this connection however is not stable under all Reidemeister
moves (e.g., one may need to create a handle to allow a
type II Reidemeister move), we will view virtual link diagrams
just abstractly as given above.

In this paper, $\oZ_+=\{0,1,2,\ldots\}$ and for any $n\in\oZ_+$:
\dez{
[n]:=\{1,\ldots,n\}.
}

Choose $n\in\oZ_+$.
Let the symmetric group $S_2$ act on $(\oR^n)^{\otimes 4}$ so that the nonidentity element of
$S_2$ brings $x_1\otimes x_2\otimes x_3\otimes x_4$
to $x_3\otimes x_4\otimes x_1\otimes x_2$.
Define
\dez{
\RR_n:=((\oR^n)^{\otimes 4})^{S_2},
}
which is the linear space of $S_2$-invariant elements of $(\oR^n)^{\otimes 4}$.
Note that $\RR_n$ can be identified with the collection of symmetric
matrices in $(\oR^{n\times n})^{\otimes 2}$.

Following de la Harpe and Jones [4], we call any element $R$ of $\RR_n$ a {\em vertex model}
(`edge-coloring model' in [12]).
For any $R\in\RR_n$, let $f_R$ be the {\em partition function} of $R$;
that is, $f_R$ is the function $f_R:\GG\to\oR$ defined by
\de{23ap07c}{
f_R(G)=\sum_{\phi:EG\to[n]}
\prod_{v\in VG}R_{\phi(\delta(v))}.
}
Here we put
\dez{
\phi(\delta(v)):=(\phi(e_1),\phi(e_2),\phi(e_3),\phi(e_4)),
}
where $e_1,e_2,e_3,e_4$ are the edges incident with $v$, in clockwise order,
and where $e_1,e_3$ form the overcrossing pair.
Since $R$ is $S_2$-invariant, $R_{\phi(\delta(v))}$ is well-defined.
Note that $f_R(\bgcirc)=n$.

The well-known sufficient conditions on $R$ for $f_R$ to be a virtual link invariant are:
\di{18fe15a}{
\nr{i} $\dps\sum_{a}R_{iaaj}=\delta_{ij}$ for all $i,j$,
\nrs{ii} $\dps\sum_{a,b}R_{ijab}R_{alkb}=\delta_{ik}\delta_{jl}$ for all $i,j,k,l$,
\nrs{iii} $\dps\sum_{a,b,c}R_{iabh}R_{jkca}R_{bclm}=\sum_{a,b,c}R_{ijbc}R_{bkla}R_{camh}$ for all $i,j,k,l,m,h$,
}
where $R$ is expressed in the standard basis of $(\oR^n)^{\otimes 4}$,
where all indices run from $1$ to $n$, and where $\delta_{ij}$ is the Kronecker delta.
Condition (iii) is the {\em Yang-Baxter equation}.
In the real case, the conditions \rf{18fe15a} are also necessary conditions for $f_R$ to be a virtual link invariant.
Elements $R$ of $\RR_n$ satisfying \rf{18fe15a} are called {\em R-matrices}.
(Often condition (i) is deleted, to obtain an invariant for `ribbon links'.)

In this paper, we characterize which real-valued functions $f$ on the collection $\GG$ are equal to $f_R$ for some
R-matrix $R$.
To this end, we introduce the concept of a $k$-join of virtual link diagrams (for any $k\in\oZ_+$).
To define it, we consider the linear space $\oR\GG$ of all formal $\oR$-linear combinations of elements
of $\GG$.
Any function on $\GG$ to a linear space can be extended uniquely to a linear function on $\oR\GG$.
The elements of $\oR\GG$ are called {\em quantum virtual link diagrams}.

The {\em $k$-join} $G\join{k}H$ of virtual link diagrams $G$ and $H$ is an element of $\oR\GG$.
It is obtained from the disjoint union
of $G$ and $H$, by taking the sum over all quantum virtual link diagrams obtained as follows:
choose distinct crossings $u_1,\ldots,u_k$ of $G$ and distinct crossings $v_1,\ldots,v_k$ of $H$,
and for each $i=1,\ldots,k$
\dy{23fe15a}{
replace
~~~
\raisebox{-.4\height}{\scalebox{0.15}{\includegraphics{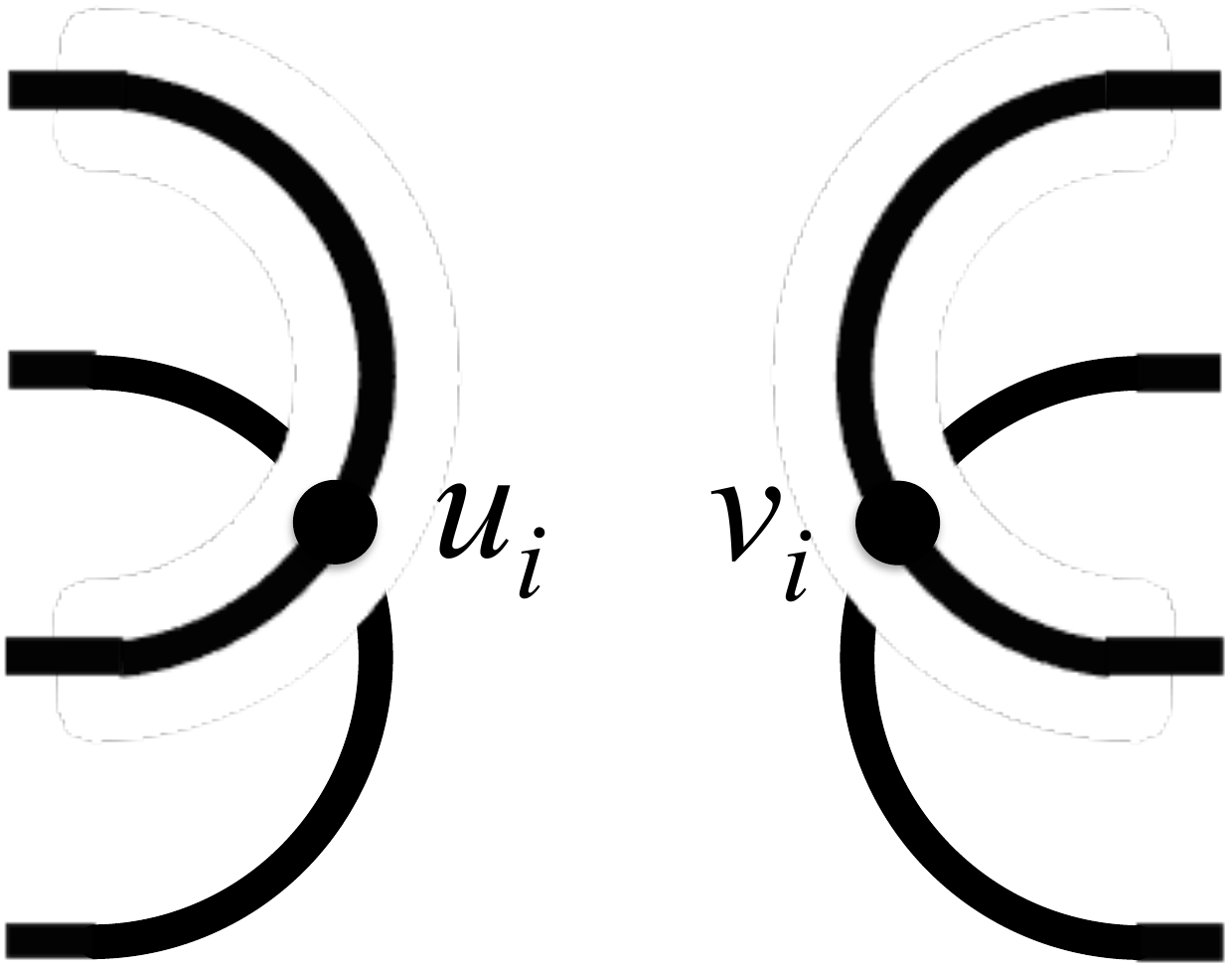}}}
~~~
by
~~~
$\frac12
\big(
~~~
\raisebox{-.42\height}{\scalebox{0.15}{\includegraphics{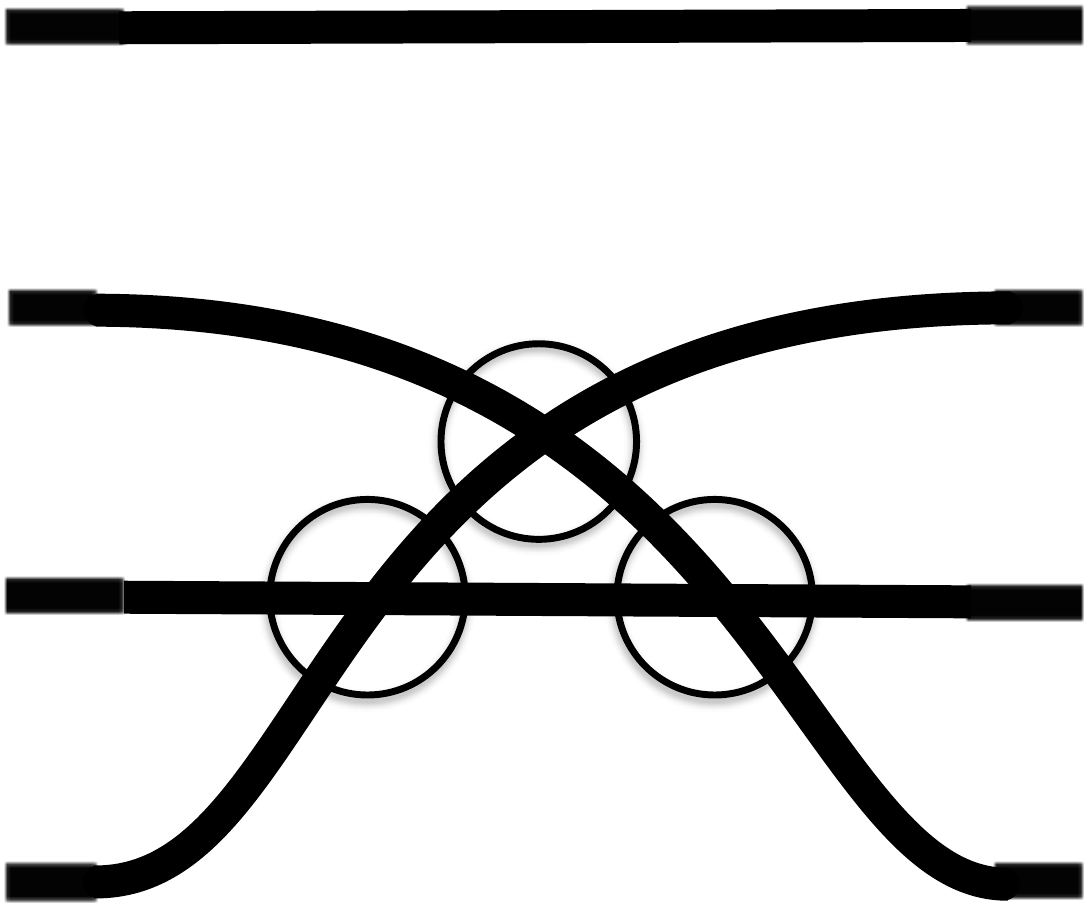}}}
~~~
+
~~~
\raisebox{-.42\height}{\scalebox{0.15}{\includegraphics{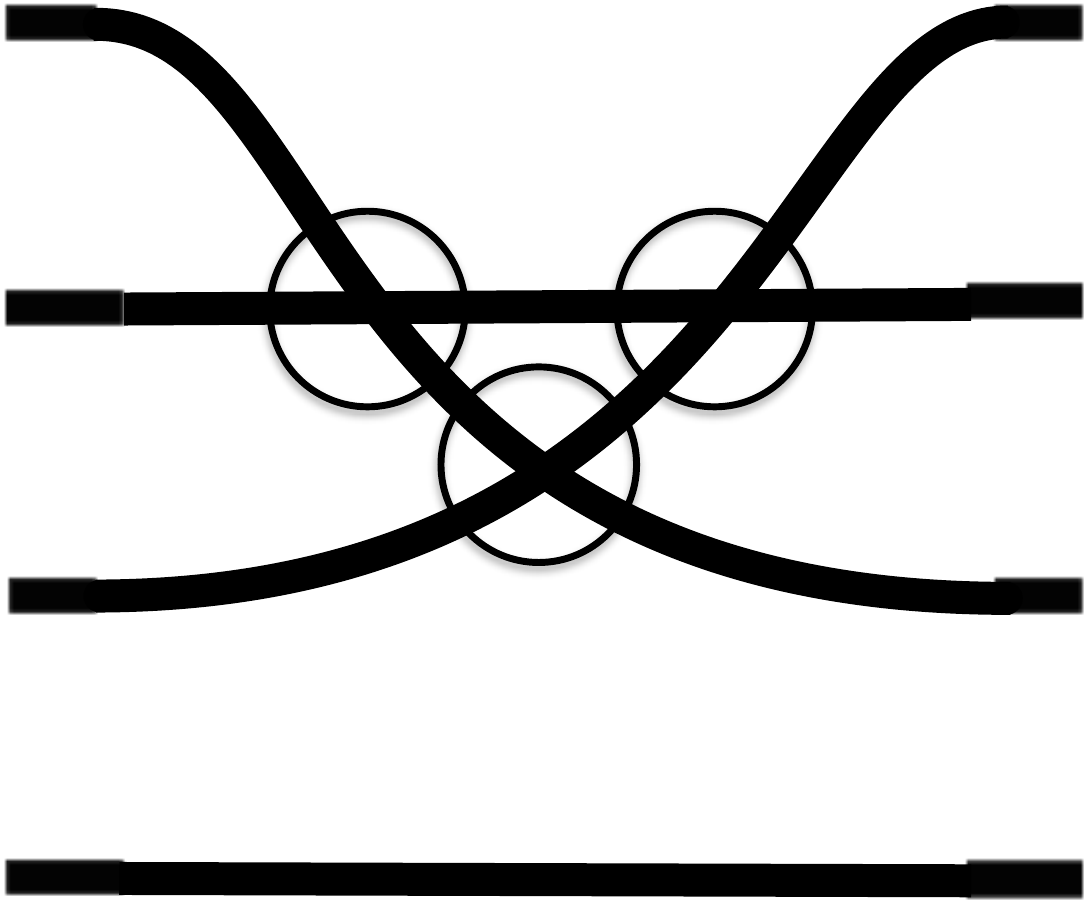}}}
~~~
\big)$.
}
As usual, a circle around a crossing in these pictures means that the crossing does not correspond to
a crossing of the virtual link diagram, but is an artefact of the planarity of the drawing.
Note that in \rf{23fe15a}, the new connections conform to the cyclic orders and the overcrossings
at $u_i$ and $v_i$.

The $k$-join can be described in terms of joining two virtual link diagram tangles (i.e., virtual link diagrams
in which labeled vertices of degree 1 are allowed) by identifying equally labeled vertices (cf.\ Szegedy [12]).
Then the $k$-join is obtained by `opening' $G$ and $H$ at the crossings $u_1,\ldots,u_k,v_1,\ldots,v_k$
(that is, deleting these vertices topologically, thus leaving, for each deleted vertex, four open end segments).
Choosing appropriate labelings at the ends and joining the tangles along equally labeled ends,
yields the $k$-join.
The $k$-join is therefore a more restricted operation, which will yield therefore a stronger characterization.

We call $f$ {\em weakly reflection positive} if for each $k\in\oZ_+$, the $\GG\times\GG$ matrix
\dyyz{
M_{f,k}:=(f(G\join{k}H))_{G,H\in\GG}
}
is positive semidefinite.
Moreover, $f:\GG\to\oR$ is called {\em multiplicative} if $f(\emptyset)=1$ (where $\emptyset$ is the virtual
link diagram with no crossings and edges) and $f(G\sqcup H)=f(G)f(H)$ for all virtual link diagrams
$G,H$, where $\sqcup$ denotes disjoint union.

\thmnn{
Let $f:\GG\to\oR$.
Then there exists an R-matrix $R$ with $f=f_R$
if and only if $f(\bgcirc)\geq 0$ and $f$ is multiplicative and weakly reflection positive and satisfies
\di{24fe15d}{
\nr{i}
$
f\big(~\raisebox{-.35\height}{\scalebox{0.10}{\includegraphics{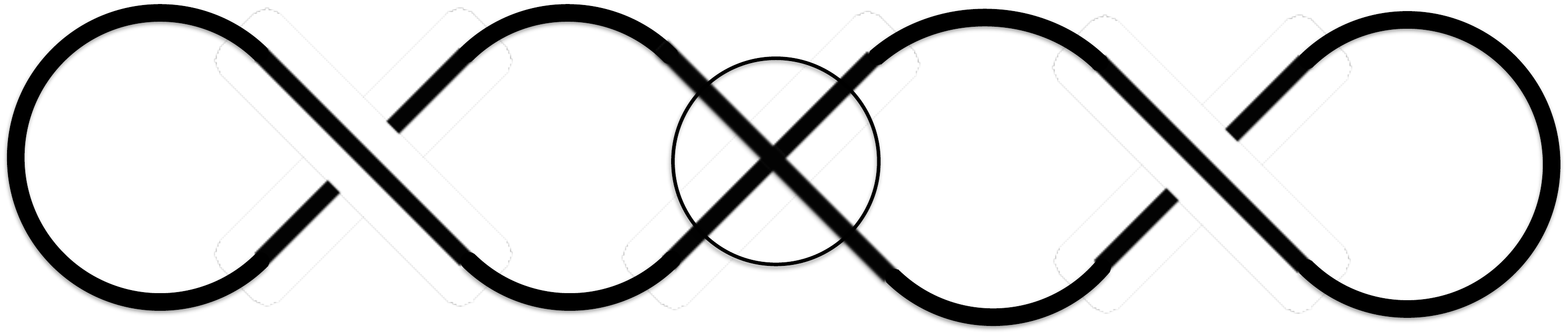}}}~\big)
+
f(\bgcirc)
=
2
f\big(~\raisebox{-.35\height}{\scalebox{0.10}{\includegraphics{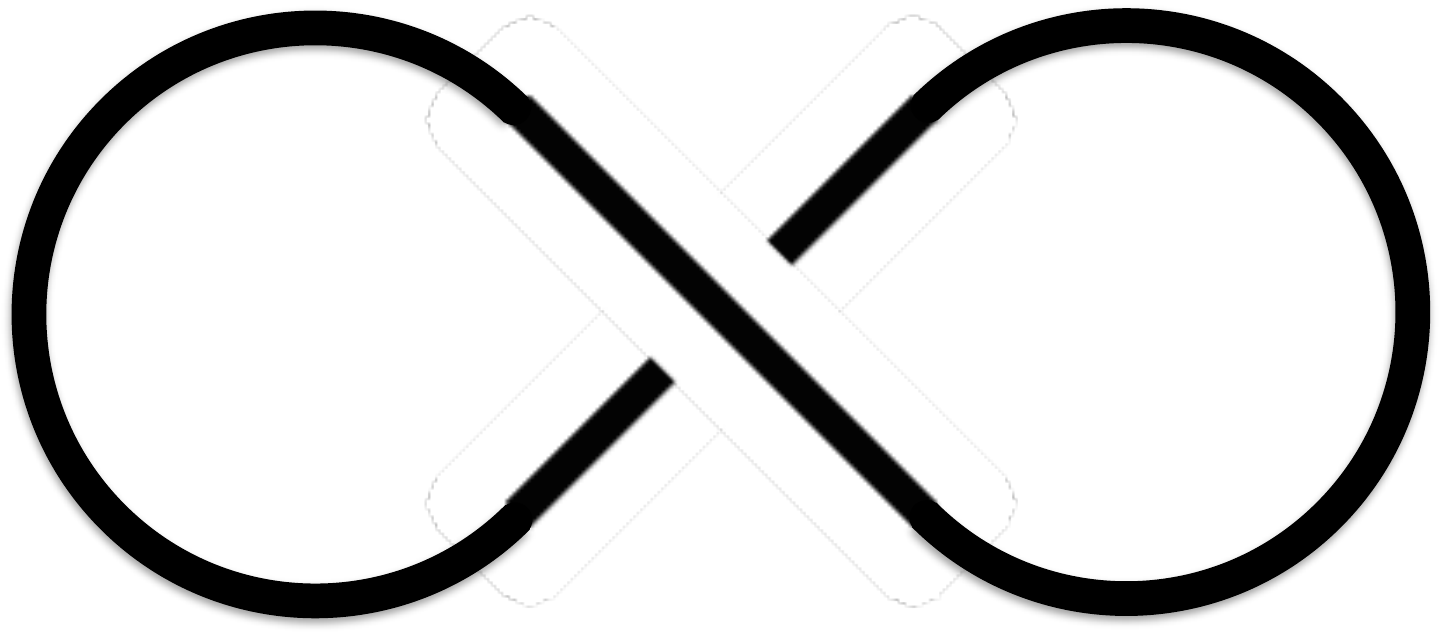}}}~\big)
$,
\nrs{ii}
$f\big(~\raisebox{-.45\height}{\scalebox{0.10}{\includegraphics{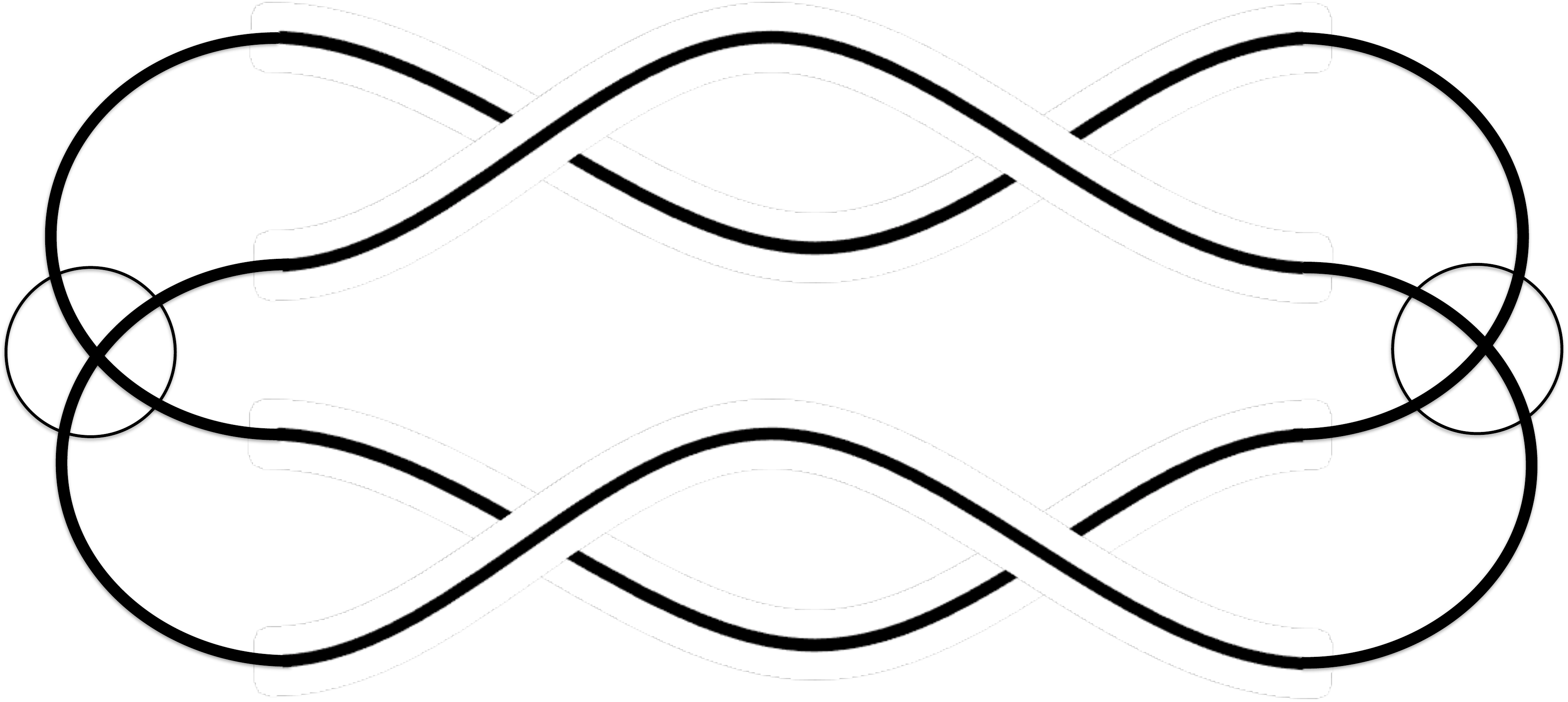}}}~\big)
+f(\bgcirc)^2
=
2f\big(~\raisebox{-.4\height}{\scalebox{0.10}{\includegraphics{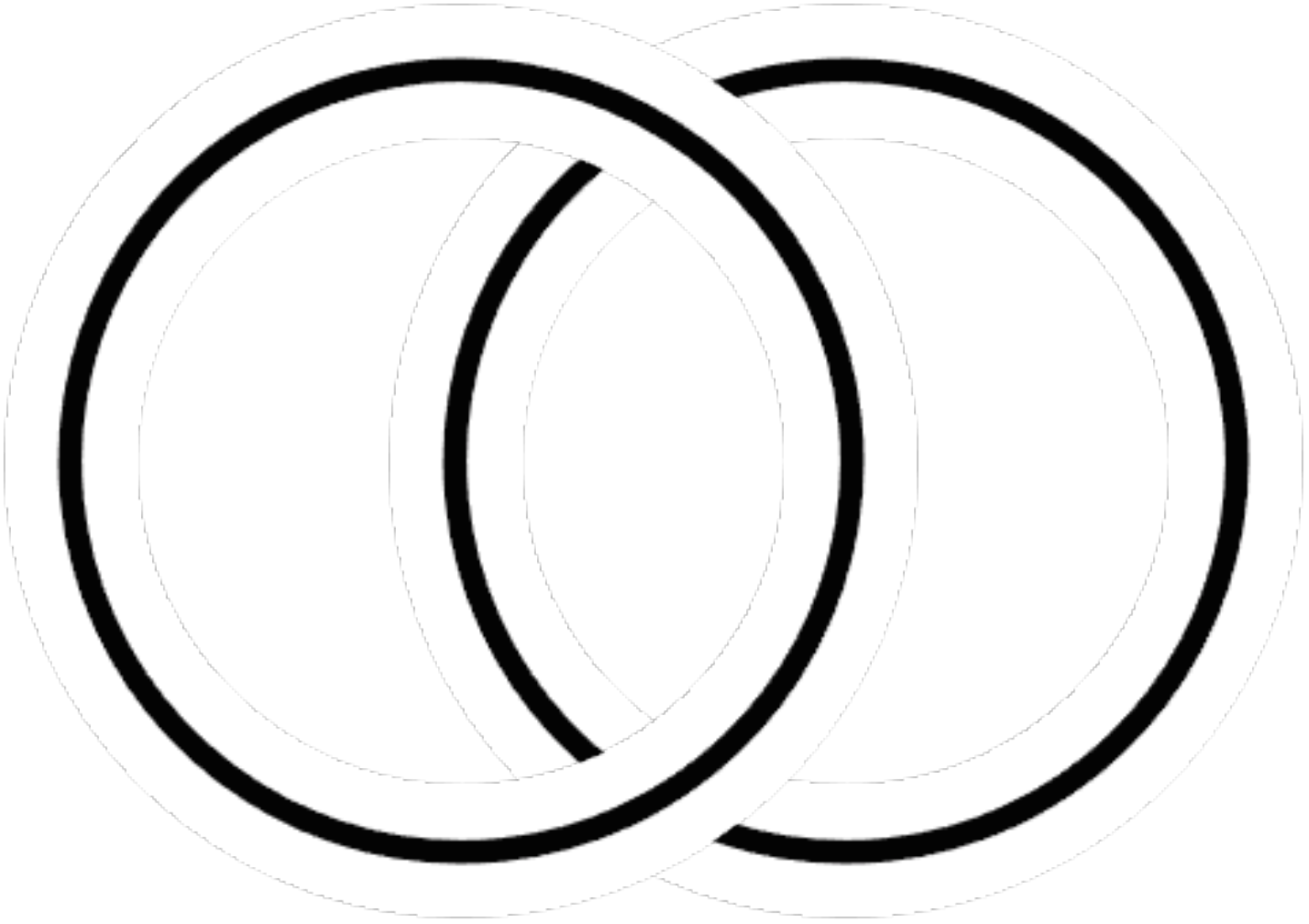}}}~\big)$,
\nrs{iii}
$f\big(~\raisebox{-.47\height}{\scalebox{0.08}{\includegraphics{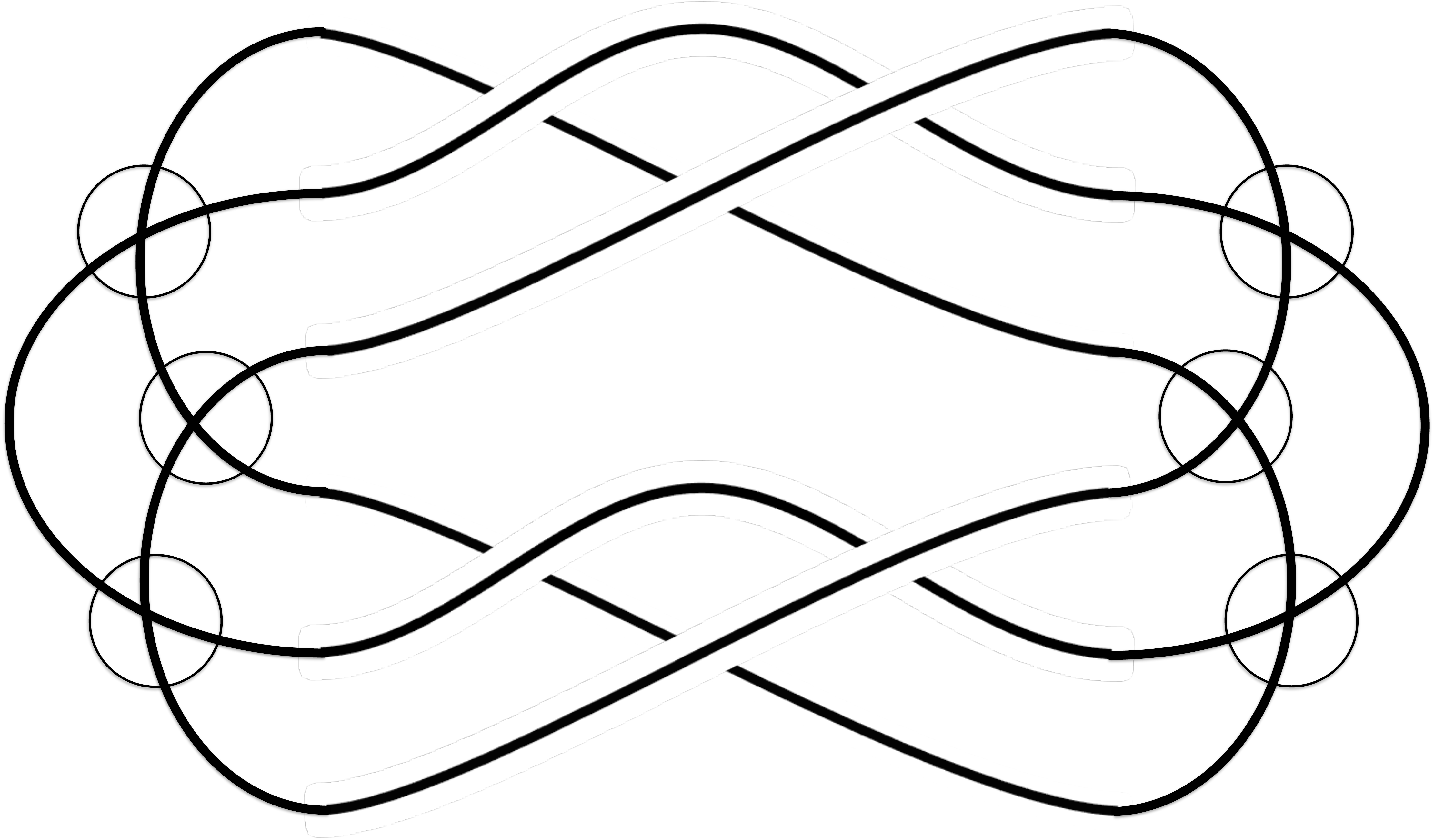}}}~\big)
=
f\big(~\raisebox{-.47\height}{\scalebox{0.08}{\includegraphics{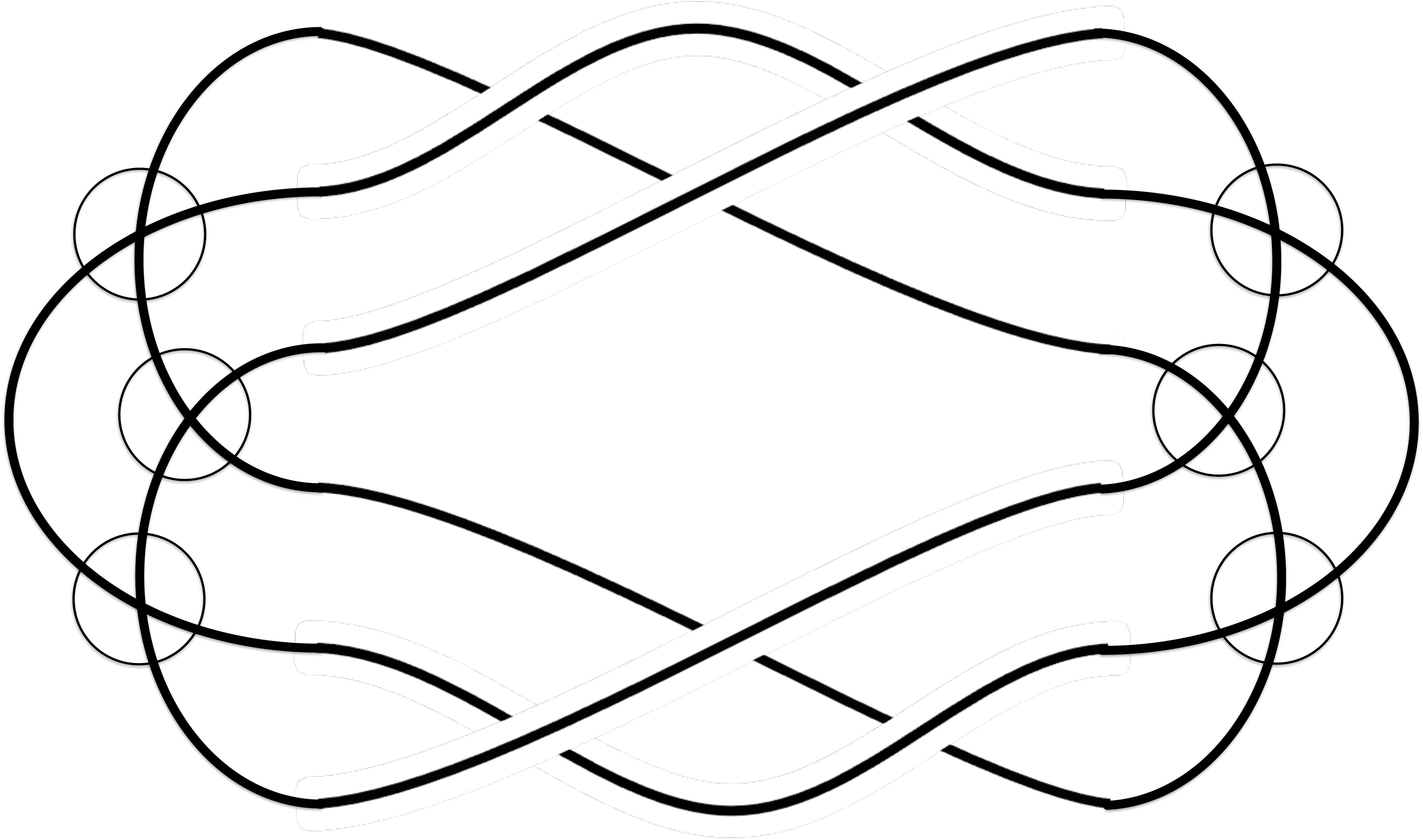}}}~\big)$.
}
}

\medskip
Our proof of the theorem follows the line of proof layed down in [9]
for `3-graphs' and cyclic cubic graphs.
The main addition of the present study is the application to virtual link diagrams, which requires a
different combinatorial proof for the integrality of $f(\bgcirc)$.
An interesting feature for virtual link diagrams is that the multiplicativity and weak reflection positivity
of $f$ imply that $f(\bgcirc)$ is an integer but might be negative.
In fact, if $f(\bgcirc)$ is negative it is even --- see the lemma below.
This raises the question to classify those multiplicative and weakly reflection positive virtual link invariants $f$
with $f(\bgcirc)<0$.

It can also be shown, with the Stone-Weierstrass theorem as in [9], that
the R-matrix $R$ in the theorem is unique, up to the natural action of the real orthogonal group $O(n)$
on $R$ (which action leaves $f_R$ invariant).

Multiplicative weakly reflection positive functions $f:\GG\to\oR$ with $f(\bgcirc)=-2k$ do exist for
any $k\in\oZ_+$.
Indeed, define $f(G)=0$ if $G$ has at least one crossing, and $f(G)=(-2k)^t$ if $G$ is the disjoint union
of $t$ copies of $\bgcirc$.
Then $f$ trivially is multiplicative, and it is weakly reflection positive, as can be derived again from the results of
Hanlon and Wales [3] displayed below.

The remainder of this paper is devoted to proving the theorem.

\sectz{The algebra homomorphism $p_n:\oR\GG\to\OO(\RR_n)$}

We make some preparations to the proof of the theorem.
The space $\oR\GG$ of formal linear combinations of elements of $\GG$, is in fact an algebra, by taking
the disjoint union $G\sqcup H$ of two virtual link diagrams $G$ and $H$ as multiplication $GH$.
Choose $n\in\oZ_+$ and recall that $\RR_n$ denotes the linear space
\dyyz{
\RR_n:=((\oR^n)^{\otimes 4})^{S_2}.
}
As usual, $\OO(\RR_n)$ denotes the algebra of polynomials on $\RR_n$.
Define an algebra homomorphism $p_n:\oR\GG\to\OO(\RR_n)$ by
\dez{
p_n(G)(R):=f_R(G)
}
for $G\in\GG$ and $R\in\RR_n$.
So the element $R$ in the theorem can be described as a common zero of the polynomials
$p_n(G)-f(G)$ for all $G\in\GG$.

We mention a connection of the $k$-join of virtual link diagrams to $k$-th derivatives of $p_n$,
which is similar to a lemma proved in [9] for cubic cyclic
graphs, and can be proved by a word for word translation of the method.

For any $q\in\OO(\RR_n)$, let $dq$ be its derivative, being an element of $\OO(\RR_n)\otimes \RR_n^*$.
So $d^kq\in\OO(\RR_n)\otimes (\RR_n^*)^{\otimes k}$.
Note that the standard inner product on $\oR^n$ induces an inner product on
$(\oR^n)^{\otimes 4}$, hence on $\RR_n$ and $\RR_n^*$, and therefore it induces a product
$\langle.,.\rangle:
(\OO(\RR_n)\otimes (\RR_n^*)^{\otimes k})\times
(\OO(\RR_n)\otimes (\RR_n^*)^{\otimes k})\to\OO(\RR_n)$.
Then, for all $G,H\in\GG$ and all $k,n\in\oZ_+$:
\de{13ap08c}{
p_n(G\join{k}H)=\langle d^kp_n(G),d^kp_n(H)\rangle.
}
This connection between $k$-joins and $k$-th derivatives will be used a number of times in our proof
of the theorem.

As in [12] (cf.\ [2],[11]), the first fundamental theorem of invariant
theory for the real orthogonal group $O(n)$ implies
\de{4ap07j}{
p_n(\oR\GG)=\OO(\RR_n)^{O(n)},
}
the latter denoting the space of $O(n)$-invariant elements of $\OO(\RR_n)$.

\sectz{The value of $f$ on $\bigcirc$}

The following lemma on $f(\bgcirc)$ carries the most combinatorial part of the proof.
It is based on basic results of Hanlon and Wales [3] on the representation theory of the symmetric group
(cf.\ Sagan [10]).

\lemmann{
If $f:\GG\to\oR$ is multiplicative and weakly reflection positive, then $f(\bgcirc)$ belongs to
$\{\ldots,-6,-4,-2,0,1,2,3,\ldots\}$.
}

\pf
I. 
We first describe some tools, using results of [3].
Consider any $k\in\oZ_+$.
For any matching $M$ on $[8k]$ and any $\pi\in S_{8k}$, let $\pi\cdot M$ be the
matching $\{\pi(e)\mid e\in M\}$.
Define $\MM$ to be the set of perfect matchings on $[8k]$. 
So the group $S_{8k}$ acts on $\MM$, which induces an action of $S_{8k}$ on $\oR^\MM$.

To each $M\in \MM$ we can associate a virtual link diagram $G_M$ on $[2k]$ by identifying,
for each $j\in[2k]$, the vertices $4j-3,4j-2,4j-1,4j$ of $M$ to one crossing called $j$
as in
\dyz{
\raisebox{-.3\height}{\scalebox{0.16}{\includegraphics{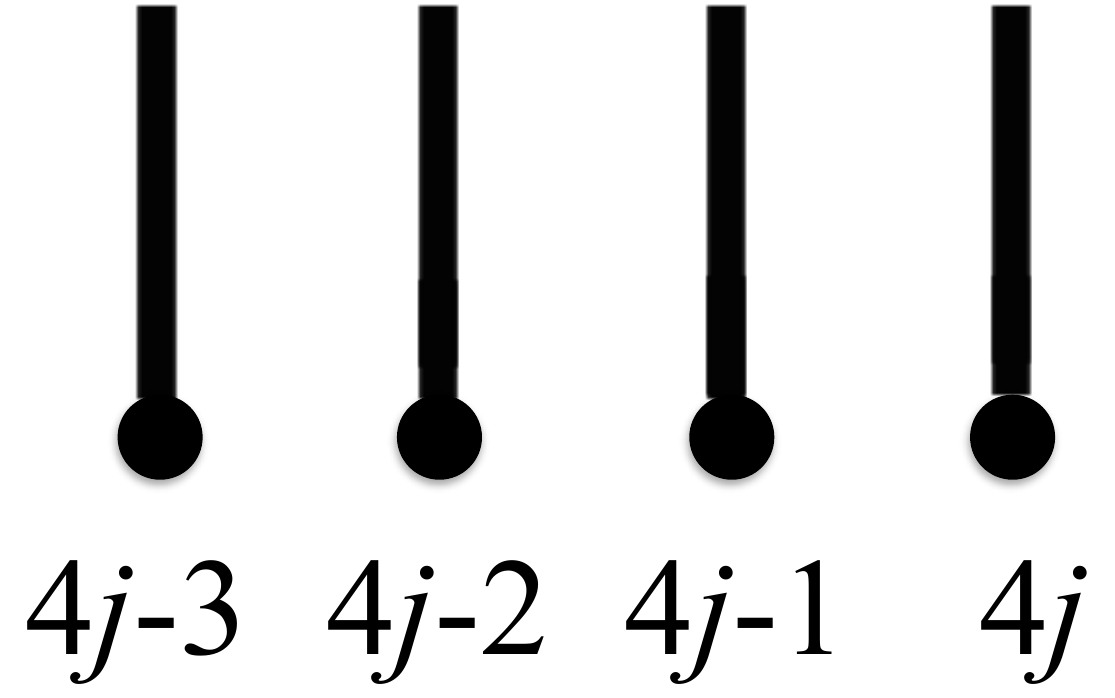}}}
~~~
$\longrightarrow$
~~~
\raisebox{-.47\height}{\scalebox{0.16}{\includegraphics{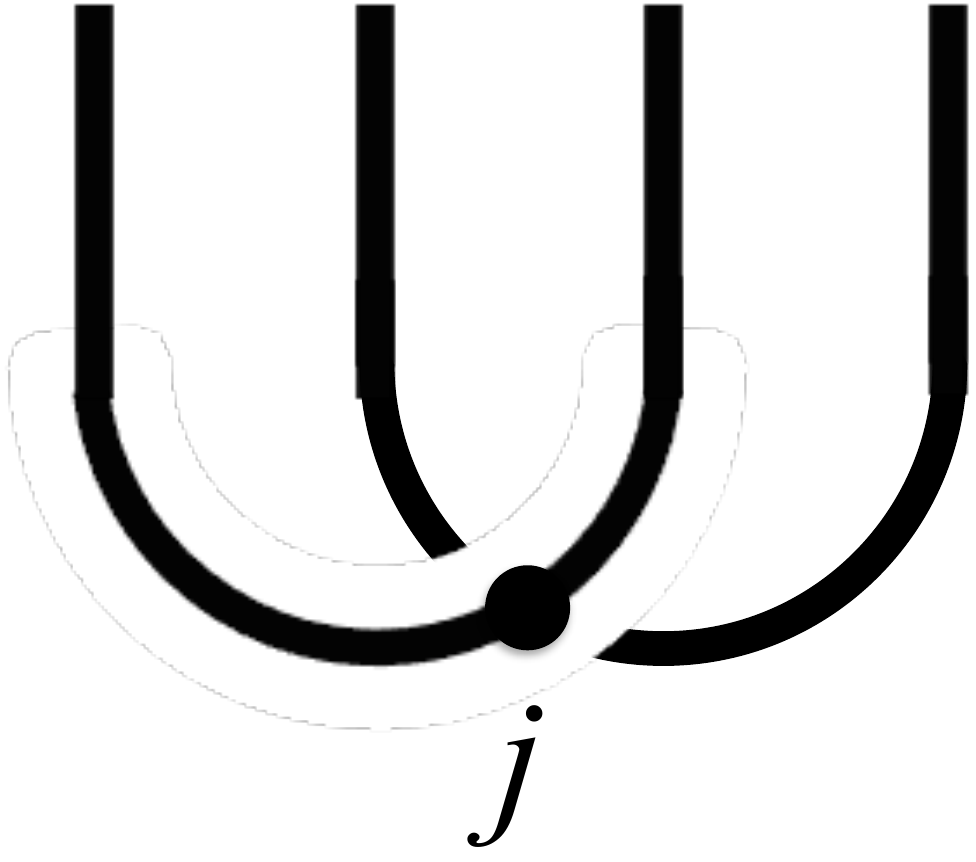}}}
~~.
}

To describe $G_M\join{2k}G_N$ for $M,N\in\MM$, we define
the following subgroups of $S_{8k}$.
For $j\in [2k]$, let $B_j$ be the group consisting of the identity $\id$ and of
$(4j-3,4j-1)(4j-2,4j)$.
Define $B := B_1B_2\cdots B_{2k}$.
Let $D$ be the group of permutations $d\in S_{8k}$ for which there exists $\pi\in S_{2k}$ such that
$d(4j-i)=4\pi(j)-i$ for each $j=1,\ldots,2k$ and $i=0,\ldots,3$.
Set $Q:=BD$, which is a group.

For $M,N\in \MM$, let $c(M,N)$ denote the number of connected components of the graph $([8k],M\cup N)$. 
Then, by definition of the operation $\join{2k}$, we have
\begin{equation}\label{eq:vee=union}
G_M \join{2k} G_N=2^{-2k}(2k)!\sum_{s\in Q}\bgcirc^{c(M,s\cdot N)}.
\end{equation}

For $\pi\in S_{8k}$, let $P_{\pi}$ be the $\MM\times\MM$ permutation matrix corresponding to $\pi$;
then $P_{\pi}w=\pi\cdot w$ for each $w\in\oR^\MM$.
For any $x\in\oR$, let $A(x)$ and $A^Q(x)$ be the $\MM\times\MM$ matrices defined by 
\dyz{
$(A(x))_{M,N} :=x^{c(M,N)}$
~~and~~
$\dps A^Q(x):=\sum_{s\in Q}A(x)P_s$,
}
for $M,N\in \MM$.
So, by the weak reflection positivity of $f$, \rf{eq:vee=union} implies that $A^Q(f(\bgcirc))$ is positive semidefinite.
Note that each $P_{\pi}$ commutes with $A(x)$, as for all $M,N\in\MM$
one has $c(\pi\cdot M,\pi\cdot N)=c(M,N)$, implying $A(x)=P_{\pi}\T A(x)P_{\pi}=P_{\pi}^{-1}A(x)P_{\pi}$.

Hanlon and Wales [3] showed that the eigenvalues and eigenvectors of $A(x)$ can be described as follows.
Consider any partition $\lambda=(t_1,\ldots,t_m)$ of $8k$, with all $t_i$ even.
Then $A(x)$ has an eigenvalue
\de{24fe15i}{
\mu_{\lambda}(x):=\prod_{a=1}^m\prod_{b=1}^{\frac12t_a}(x-a+2b-1).
}
To describe a corresponding eigenvector, make a Young tableau $T$ associated to
$\lambda$ such that each row of $T$ has the form
\renewcommand{\arraystretch}{1.2}
\dez{
\begin{array}{|c|c|c|c|c|c|c|}
\hline
i_1&\overline{i_1}&
i_2&\overline{i_2}&
\cdots&
i_t&\overline{i_t}\\
\hline
\end{array}
}
for some $i_1,\ldots,i_t\in[4k]$, where $\overline{i}:=4k+i$ for each $i\in[4k]$.
For $i=1,\ldots,t_1$, let $K_i$ denote the set of numbers in column $i$ of $T$ and let $C_i$ be
the subgroup of $S_{8k}$ that permutes the elements of $K_i$.
Then $C:=C_1\cdots C_{t_1}$.
Similarly, for $i=1,\ldots,m$, let $R_i$ be the subgroup of $S_{8k}$ that permutes the numbers in row $i$ of $T$,
and $R:=R_1\ldots R_m$.

Let $F$ be the perfect matching on $[8k]$ with edges $\{i,\overline i\}$ for $i\in[4k]$.
Then
\dez{
v:=\sum_{c\in C,r\in R}\sgn(c)cr\cdot F
}
is an eigenvector of $A(x)$ belonging to $\mu_{\lambda}(x)$.
Then for $u:=\sum_{q\in Q}q\cdot v$ one has
\dyyz{
\hspace*{-15pt}A^Q(x)u=\sum_{q',q\in Q}A(x)P_{q'}P_{q}v=\sum_{q',q\in Q}P_{q'}P_{q}A(x)v=\mu_{\lambda}(x)\sum_{q',q\in Q}P_{q'}P_{q}v=|Q|\mu_{\lambda}(x) u.
}
So $u$ is an eigenvector of $A^Q(x)$ belonging to $|Q|\mu_{\lambda}(x)$, {\em provided that} $u$ is
nonzero.
For this it suffices that the coefficient $u_F$ of $u$ in $F$ is nonzero.
Note that
\dyyz{
u_F=\sum_{q\in Q}(q\cdot v)_F=\sum_{q\in Q}\sum_{c\in C,r\in R}\sgn(c)(qcr\cdot F)_F
=
\sum_{q\in Q,c\in C,r\in R\atop qcr\cdot F=F}\sgn(c).
}
So $u\neq 0$ if
for any $q\in Q$, $c\in C$, and $r\in R$, if $qcr\cdot F=F$ then $\sgn(c)=1$;
that is (as $Q$ is a group), if for any $q\in Q$, $c\in C$, $r\in R$:
\dy{24fe15h}{
if $q\cdot F=cr\cdot F$, then $\sgn(c)=1$.
}

\medskip
\noindent
II. We first apply part I to the case where $f(\bgcirc)\geq 0$.
Let $k:=\lceil f(\bgcirc)\rceil+1$, and consider the partition $\lambda:=(8,8,\ldots,8)$ of $8k$.
Then, by \rf{24fe15i},
\dyyz{
\mu_{\lambda}(x)=\prod_{i=0}^{k-1}(x-i)(x-i+2)(x-i+4)(x-i+6).
}
We give a Young tableau associated to $\lambda$ that will yield \rf{24fe15h}.
This implies that $|Q|\mu_{\lambda}(x)$ is an eigenvalue of $A^Q(x)$.
So $\mu_{\lambda}(f(\bgcirc))\geq 0$.
Hence, as the polynomial $\mu_{\lambda}(x)$ has largest zero $k-1$, with multiplicity 1, and
as $k-1=\lceil f(\bgcirc)\rceil$, we know $f(\bgcirc)=k-1$.

Consider the following Young tableau associated to $\lambda$:
\dyyz{
T:=\begin{array}{|c|c|c|c|c|c|c|c|}
\hline
1&\overline 1&2&\overline 2&3&\overline 3&4&\overline 4\\
\hline
5&\overline 5&6&\overline 6&7&\overline 7&8&\overline 8\\
\hline
\vdots&\vdots&\vdots&\vdots&\vdots&\vdots&\vdots&\vdots\\
\hline
4k\minus 3&\overline{4k\minus 3}& 4k\minus 2&\overline{4k\minus 2}& 4k\minus 1&\overline{4k\minus 1}& 4k&\overline{4k}\\
\hline
\end{array}.
}
To prove \rf{24fe15h}, choose $q\in Q$, $c\in C$, and $r\in R$ with $q\cdot F=cr\cdot F$.
Let $c=c_1\cdots c_8$ with $c_i\in C_i$ ($i=1,\ldots,8$) and define $M:=q\cdot F$.
Since $F$ has no edges between $X:=K_1\cup K_2\cup K_5\cup K_6$ (the set of odd numbers in $T$)
and $Y:=K_3\cup K_4\cup K_7\cup K_8$ (the set of even numbers in $T$) and
since $Q\cdot X=X$ and $Q\cdot Y=Y$, we know that $M$ has no edges between $X$ and $Y$.
For any $N\in\MM$ and $Z\subseteq [8k]$, let $N_Z$ be the set of edges of $N$ contained in $Z$.

Let $z\in S_{8k}$ be defined by $z(i):=i+1$ if $4$ does not divide $i$ and $z(i):=i-3$
if 4 divides $i$.
So $z^4=\id$, $z(X)=Y$, and $z\cdot F=F$.
Moreover, $zq=qz$ (since $zb=bz$ and $zd=dz$ for all $b\in B$ and $d\in D$).
So $z\cdot M=M$.
Hence $z\cdot M_X=M_Y$.

Let $N:=r\cdot F$.
So $M=c\cdot N$.
As no edge of $M$ connects $X$ and $Y$, also no edge in $N$ connects $X$ and $Y$.
Moreover, as $z\cdot M_X=M_Y$, for each two columns $K_i$ and $K_j$ in $X$, we have
$|M_{K_i\cup K_j}|=|M_{K_{i+2}\cup K_{j+2}}|$, and hence
$|N_{K_i\cup K_j}|=|N_{K_{i+2}\cup K_{j+2}}|$.
Moreover, if an edge $e\in N$ connects $K_i$ and $K_j$, then $N$ has an edge in the same
row as $e$ connecting the other two columns in $X$; similarly for $Y$.

This implies that there exists a permutation $c'\in C_1C_2C_5C_6$ that permutes complete rows in $X$
in such a way that $c'\cdot N_X$ is a shift of $N_Y$; that is, $zc'\cdot N_X=N_Y$.
As $c'$ maintains rows in $X$, there exists $r'\in R$ with $c'\cdot N=r'\cdot F$; so $c(c')^{-1}r'\cdot F=cr\cdot F$.
Moreover, $\sgn(c')=1$,
and, setting $N':=r'\cdot F$ we have
$z\cdot N'_X=z\cdot(r'\cdot F)_X=z\cdot(c'\cdot N)_X=zc'\cdot N_X=N_Y=N'_Y$.
Therefore, by replacing $r$ by $r'$ and $c$ by $c(c')^{-1}$ we can assume that $z\cdot N_X=N_Y$.

Next consider any two columns $K_i$ and $K_j$ in $X$.
Let $X':=K_i\cup K_j$ and $Y':=K_{i+2}\cup K_{j+2}$. So $Y'=z(X')$ and $z\cdot N_{X'}=N_{Y'}$.
Then $e\mapsto z^{-1}c^{-1}zc(e)$ is a permutation $\sigma$ of the edges $e$ in $N_{X'}$, since
$e\in N_{X'}$ $\Rightarrow$ $c(e)\in M_{X'}$ $\Rightarrow$ $zc(e)\in M_{Y'}$ $\Rightarrow$ $c^{-1}zc(e)\in N_{Y'}$ $\Rightarrow$ $z^{-1}c^{-1}zc(e)\in z^{-1}\cdot N_{Y'}=N_{X'}$.
As $\sigma$ permutes edges in $X'$, there exists a permutation $c'\in C_iC_j$ such that
$c'(e)=z^{-1}c^{-1}zc(e)$ for all $e\in N_{X'}$ and such that $c'$ only permutes elements covered by $N_{X'}$.
Then $\sgn(c')=1$.
By replacing $c$ by $c(c')^{-1}$ we attain that $e=z^{-1}c^{-1}zc(e)$ for all edges $e\in N_{X'}$.
So $cz(e)=zc(e)$ for all $e\in N_{X'}$.

Doing this for all $K_i$ and $K_j$ in $X$, we finally achieve that $cz(e)=zc(e)$ for all $e\in N_X$.
As $N_X$ is a perfect matching on $X$, this implies $cz(i)=zc(i)$ for all $i\in X$.
Equivalently, $c_3c_4c_7c_8z(i)=zc_1c_2c_5c_6(i)$ for all $i\in X$.
Hence $\sgn(c_3c_4c_7c_8)=\sgn(c_1c_2c_5c_6)$, implying $\sgn(c)=1$.

\medskip
\noindent
III. Next we apply part I of this proof to the case where $f(\bgcirc)\leq 0$.
Choose $k\in\oZ_+$, and consider the partition $\lambda:=(8k)$ of $8k$ and the following Young tableau
\dyyz{
T:=\begin{array}{|c|c|c|c|c|c|c|c|c|c|c|c|c|c|c|c|c|c|c|}
\hline
1&\overline 1&
2&\overline 2&
\cdots&
4k\minus 1&\overline{4k\minus 1}&
4k&\overline{4k}\\
\hline
\end{array}~.
}
Then by \rf{24fe15i},
\dez{
\mu_{\lambda}(x)=\prod_{b=1}^{4k}(x-2+2b).
}
Moreover, \rf{24fe15h} trivially holds, as $C$ only consists of the identity.
The zeros of $\mu_{\lambda}$ are $-8k+2,-8k+4,-8k+6,\ldots,-2,0$, all with multiplicity 1, so that
$\mu_{\lambda}(f(\bgcirc))\geq 0$ implies that $f(\bgcirc)$ does not belong to any interval
$(-4t-2,-4t)$ for any $t\in\oZ_+$ with $t<2k$.
As $k$ can be chosen arbitrarily large, we know that $f(\bgcirc)\not\in(-4t-2,-4t)$ for all $t\in\oZ_+$.

To exclude the intervals $(-4t-4,-4t-2)$, consider the partition $\lambda:=(8k-2,2)$ of $8k$ and the
Young tableau
\dyyz{
T:=\begin{array}{|c|c|c|c|c|c|c|c|c|c|c|c|c|c|c|c|c|}
\hline
1&\overline 1&
3&\overline 3&
4&\overline 4&
\cdots&
4k\minus 1&\overline{4k\minus 1}&
4k&\overline{4k}\\
\hline
2&\overline 2\\
\cline{1-2}
\end{array}~.
}
In this case, by \rf{24fe15i},
\dez{
\mu_{\lambda}(x)=(x-1)\prod_{b=1}^{4k-1}(x-2+2b).
}
To show \rf{24fe15h}, let $c=c_1c_2$ with $c_1\in C_1$, $c_2\in C_2$.
Observe that $M:=q\cdot F$ contains no edges connecting an odd number with an even number
(as $F$ does not, and as $Q$ maintains the sets of odd and even numbers).

If $\{2,\overline 2\}$ belongs to $M$, then either $c_1$ and $c_2$ both are the identity permutation,
or $c_1$ and $c_2$ both are transpositions.
In either case, $\sgn(c)=1$ follows.

If $\{2,\overline 2\}$ does not belong to $M$, then $2$ and $\overline 2$ are matched in $M$ to even
numbers in the first row of $T$.
In this case, both $c_1$ and $c_2$ are transpositions, and again $\sgn(c)=1$ follows.
This proves \rf{24fe15h}.

Now the zeros of $\mu_{\lambda}$ are $-8k+4,-8k+6,\ldots,-2,0,1$, all with multiplicity 1,
so that, like above, $f(\bgcirc)\not\in(-4t-4,-4t-2)$ for all $t\in\oZ_+$.
\bx

\sectz{Proof of the theorem}

To see necessity in the theorem, let $R$ be an R-matrix, say $R\in\RR_n$.
Then $f_R$ is trivially multiplicative.
Positive semidefiniteness of $M_{f_R,k}$ follows from
\dyyz{
f_R(G\join{k}H)=
p_n(G\join{k}H)(R)=
\langle d^kp_n(G)(R),d^kp_n(H)(R)\rangle,
}
using \rf{13ap08c}.

To prove sufficiency, let $f$ satisfy the conditions of the theorem.
As $f(\bgcirc)\geq 0$ by assumption, the lemma implies that $n:=f(\bgcirc)$ is a nonnegative integer.
Then
\dy{20ap08b}{
there exists an algebra homomorphism $F:p_n(\oR\GG)\to\oR$ such that $f=F\circ p_n$.
}
Otherwise, as $f$ and $p_n$ are algebra homomorphisms, there exists a quantum virtual link diagram
$\gamma$ with $p_n(\gamma)=0$
and $f(\gamma)\neq 0$.
We can assume that $p_n(\gamma)$ is homogeneous, that is, all virtual link diagrams
in $\gamma$ have the same number of crossings, $k$ say.
So $\gamma\join{k}\gamma$ has no crossings, that is, it is a polynomial in $\bgcirc$.
As moreover $f(\bgcirc)=n=p_n(\bgcirc)$,
we have $f(\gamma\join{k}\gamma)=p_n(\gamma\join{k}\gamma)=0$,
the latter equality because of \rf{13ap08c}.
Similarly to Lemma 1 of [9], $\gamma$ belongs to
the ideal in $\oR\GG$ generated by $\gamma\join{k}\beta^i$ ($i=0,\ldots,k$), where
$\beta$ is the virtual link diagram
\dyy{24fe15b}{
\beta:=
~
\raisebox{0.9cm}{\rotatebox{270}{\scalebox{0.10}{\includegraphics{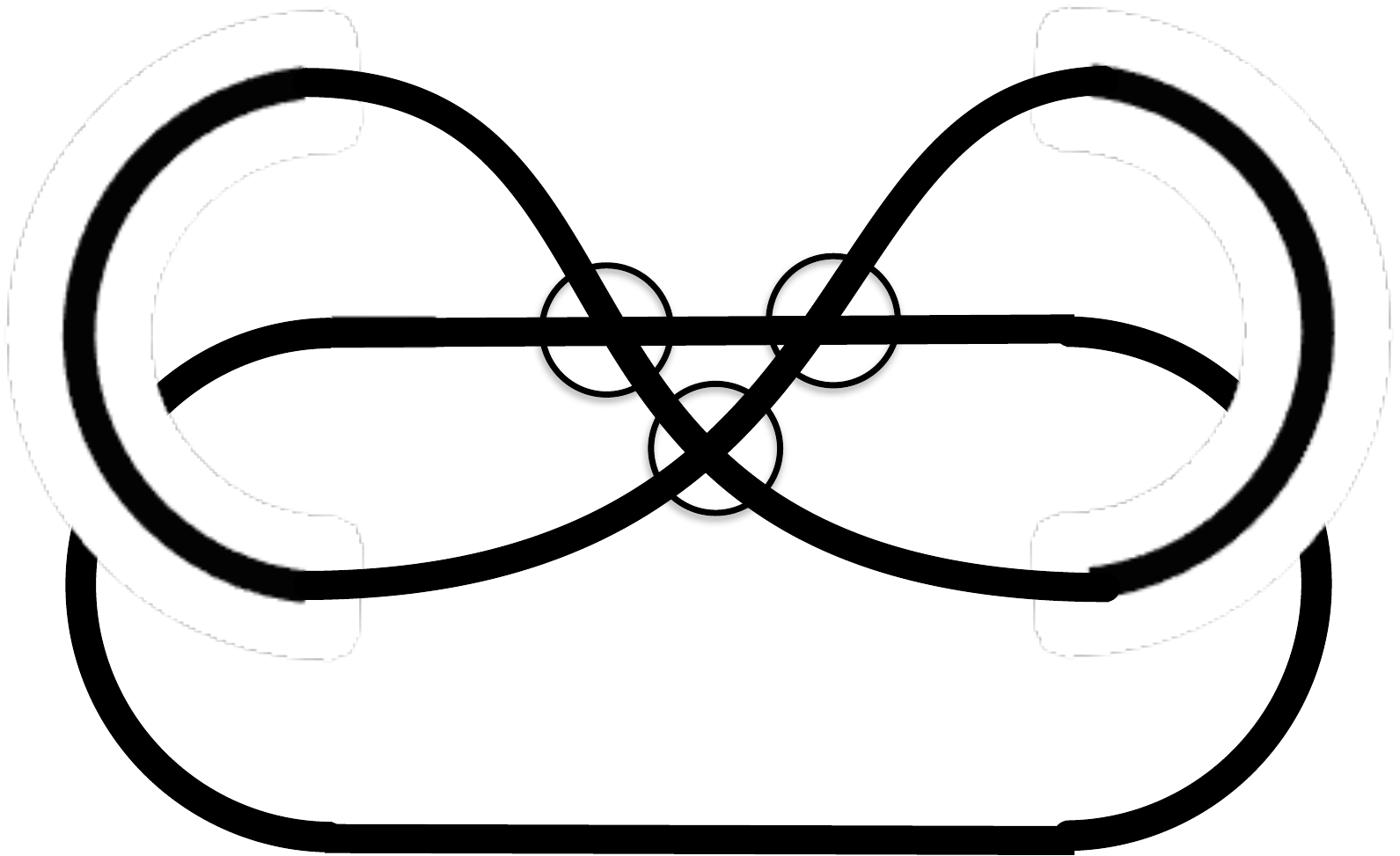}}}}~~.
}
(Note that $G\join{1}\beta=2|V(G)|G$ for each virtual link diagram $G$.)
As $f(\gamma\join{k}\gamma)=0$ implies that $f(\gamma\join{k}\beta^i)=0$ for each $i$ (by the weak reflection
positivity of $f$), we know $f(\gamma)=0$, proving \rf{20ap08b}.

Now, by \rf{4ap07j},
$p_n(\oR\GG)=\OO(\RR_n)^{O(n)}$.
Basic invariant theory then gives the existence of an $R$ in the complex extension of $\RR_n$ such that
$F(q)=q(R)$ for each $q\in\OO(\RR_n)^{O(n)}$ (cf.\ [9]).
To prove that we can take $R$ real, we apply the Procesi-Schwarz theorem [8].

For all $G,H\in\GG$, using \rf{13ap08c}:
\de{17ap08b}{
F(\langle dp_n(G),dp_n(H)\rangle)
=
F(p_n(G\join{1}H))
=
f(G\join{1}H)
=
(M_{f,1})_{G,H}.
}
Since $M_{f,1}$ is positive semidefinite, \rf{17ap08b} implies
$F(\langle dq,dq\rangle)\geq 0$ for each $q\in p_n(\oR\GG)=\OO(\RR_n)^{O(n)}$.
Then by [8] there exists a (real) $R\in \RR_n$
such that $F(q)=q(R)$ for each $q\in\OO(\RR_n)^{O(n)}=p_n(\oR\GG)$.
Then $f=f_R$, as $f(G)=F(p_n(G))=p_n(G)(R)=f_R(G)$ for each $G\in\GG$.

One may finally check that substituting $f:=f_R$ in \rf{24fe15d}, condition \rf{24fe15d}(i) is equivalent to
\dez{
\sum_{i,j}\big(\sum_{a}R_{iaaj}-\delta_{ij}\big)^2=0,
}
and hence to \rf{18fe15a}(i);
condition \rf{24fe15d}(ii) is equivalent to
\dez{
\sum_{i,j,k,l}\big(\sum_{a,b}R_{ijab}R_{alkb}-\delta_{ik}\delta_{jl}\big)^2=0,
}
and hence to \rf{18fe15a}(ii);
and condition \rf{24fe15d}(iii) is equivalent to
\dez{
\sum_{i,j,k,l,m,h}\big(\sum_{a,b,c}R_{iabh}R_{jkca}R_{bclm}-\sum_{a,b,c}R_{ijbc}R_{bkla}R_{camh}\big)^2=0,
}
and hence to \rf{18fe15a}(iii).
So $R$ is an R-matrix, as required.
\bx

\section*{References}\label{REF}
{\small
\begin{itemize}{}{
\setlength{\labelwidth}{4mm}
\setlength{\parsep}{0mm}
\setlength{\itemsep}{0mm}
\setlength{\leftmargin}{5mm}
\setlength{\labelsep}{1mm}
}
\item[\mbox{\rm[1]}] M.H. Freedman, L. Lov\'asz, A. Schrijver, 
Reflection positivity, rank connectivity, and homomorphisms of graphs,
{\em Journal of the American Mathematical Society} 20 (2007) 37--51.

\item[\mbox{\rm[2]}] R. Goodman, N.R. Wallach, 
{\em Symmetry, Representations, and Invariants},
Springer, Dordrecht, 2009.

\item[\mbox{\rm[3]}] P. Hanlon, D. Wales, 
On the decomposition of Brauer's centralizer algebras,
{\em Journal of Algebra} 121 (1989) 409--445.

\item[\mbox{\rm[4]}] P. de la Harpe, V.F.R. Jones, 
Graph invariants related to statistical mechanical models:
examples and problems,
{\em Journal of Combinatorial Theory, Series B} 57 (1993) 207--227.

\item[\mbox{\rm[5]}] L.H. Kauffman, 
Virtual knot theory,
{\em European Journal of Combinatorics} 20 (1999) 663--690.

\item[\mbox{\rm[6]}] L.H. Kauffman, 
Introduction to virtual knot theory,
{\em Journal of Knot Theory and Its Ramifications} 21 (2012) 1240007 (37 pp).

\item[\mbox{\rm[7]}] V.O. Manturov, D.P. Ilyutko, 
{\em Virtual Knots --- The State of the Art},
World Scientific, River Edge, N.J., 2013.

\item[\mbox{\rm[8]}] C. Procesi, G. Schwarz, 
Inequalities defining orbit spaces,
{\em Inventiones Mathematicae} 81 (1985) 539--554.

\item[\mbox{\rm[9]}] G. Regts, A. Schrijver, B. Sevenster, 
On partition functions for 3-graphs,
preprint, 2015, ArXiv 1503.00337v1

\item[\mbox{\rm[10]}] B.E. Sagan, 
{\em The Symmetric Group: Representations, Combinatorial Algorithms, and Symmetric Functions},
Graduate Texts in Mathematics, Vol. 203, Springer, New York, 2001.

\item[\mbox{\rm[11]}] A. Schrijver, 
On virtual link invariants,
2012, ArXiv 1211.3572

\item[\mbox{\rm[12]}] B. Szegedy, 
Edge coloring models and reflection positivity,
{\em Journal of the American Mathematical Society}
20 (2007) 969--988.

\end{itemize}
}

\end{document}